\documentclass[12pt]{amsart}

\usepackage{amsmath,amsthm,amssymb,amsfonts}
\usepackage{mathrsfs}

\begin{document}

\title[Marginal Likelihood]{On the computation of the marginal likelihood}

\author[Marques]{Paulo C. Marques F.}
\address{Instituto de Matem\'atica e Estat\'istica da Universidade de S\~ao Paulo}
\email{pmarques@ime.usp.br}
\date{April 28, 2013}

\begin{abstract}
We describe briefly in this note a procedure for consistently estimating the marginal likelihood of a statistical model through a sample from the posterior distribution of the model parameters.
\end{abstract}

\maketitle

\thispagestyle{empty}

Consider a Bayesian analysis in which observables $X=(X_1,\dots,X_n)$ are modeled conditionally, given that parameter $\Theta=\theta$, with density $f_{X\mid\Theta}(\;\cdot\mid \theta)$. Let $\pi(\,\cdot\,)$ denote the prior density of $\Theta$. For data $x=(x_1,\dots,x_n)$, the marginal likelihood
$$
  f_X(x) = \int f_{X\mid\Theta}(x\mid \theta)\,\pi(\theta) \, d\theta
$$
is a quantity commonly used in Bayesian model selection (see ref. \cite{ohagan}) to assess the relative adequacy of the model with respect to a set of competing models of the same experiment.

If $\pi(\;\cdot\mid x)>0$ is the posterior density of the parameter $\Theta$ derived from Bayes's theorem, rewriting the marginal likelihood as
$$
  f_X(x) = \int \left(\frac{f_{X\mid\Theta}(x\mid \theta)\,\pi(\theta)}{\pi(\theta\mid x)}\right) \pi(\theta\mid x)\,d\theta \, ,
$$
we can reinterpret it as the conditional expectation
$$
  f_X(x) = \mathrm{E}\!\left[ \frac{f_{X\mid\Theta}(x\mid\Theta)\,\pi(\Theta)}{\pi(\Theta\mid x)}\;\Bigg\vert\; X=x\right] \, .
$$

The difficulty with this conditional expectation is that, in general, we don't know the expression of the posterior density $\pi(\;\cdot\mid x)$. But suppose that we have available an iid sample $\theta^{(1)},\dots,\theta^{(N)}$ from the posterior distribution of $\Theta$, obtained by any simulation method. With the help (see ref. \cite{silverman}) of a pointwise consistent kernel density estimator $\hat{\pi}_N(\;\cdot\mid x)$ of the posterior density, we get, by the Law of Large Numbers\footnote{In a non-iid MCMC setting, strong convergence is guaranteed by the ergodic theorem.}, a consistent estimator of the marginal likelihood. For the sample $\theta^{(1)},\dots,\theta^{(N)}$, the estimate is
$$
  \widehat{f_X(x)} = \frac{1}{N}\sum_{i=1}^N \frac{f_{X\mid\Theta}\left(x\;\Big\vert\; \theta^{(i)}\right) \,\pi\left(\theta^{(i)}\right)}{\hat{\pi}_N\left(\theta^{(i)}\;\Big\vert\; x\right)} \, .
$$

As a simple example, let $X_1,\dots,X_n$ be conditionally independent and identically distributed, given that $\Theta=\theta$, such that $X_1$ has distribution $\mathrm{N}(\theta,\sigma^2)$. Suppose that the variance $\sigma^2$ is known. Then,
$$
  f_{X\mid\Theta}(x\mid\theta) = (2\pi)^{-n/2} \sigma^{-n} \exp \left( -\frac{n}{2\sigma^2} \left((\theta-\bar{x})^2+s^2 \right) \right) \, ,
$$
in which $\bar{x}=\sum_{i=1}^n x_i/n$, and $s^2 = \sum_{i=1}^n (x_i-\bar{x})/n$. If $\Theta$ has prior distribution $\mathrm{N}(\theta_0,\sigma_0^2)$, a straightforward calculation using Bayes's theorem shows that the posterior distribution of $\Theta$ is
$$
\mathrm{N} \left( \frac{n\sigma_0^2\bar{x}+\sigma^2\theta_0}{n\sigma_0^2+\sigma^2}, \frac{\sigma^2\sigma_0^2}{n\sigma_0^2+\sigma^2} \right) \, ,
$$
and the marginal likelihood is
\begin{align*}
  f_X(x) =& \frac{(2\pi)^{-n/2} \sigma^{-n+1}}{\sqrt{n\sigma_0^2+\sigma^2}} \\
  & \times \exp \left( -\frac{1}{2} \left(\frac{ns^2}{\sigma^2} - \frac{(n\sigma_0^2\bar{x}+\sigma^2\theta_0)^2}{\sigma^2\sigma_0^2(n\sigma_0^2+\sigma^2)} + \frac{n\sigma_0^2\bar{x}^2+\sigma^2\theta_0^2}{\sigma^2\sigma_0^2} \right) \right)\, .
\end{align*}

Written in the free R programming language (see ref. \cite{rproject}), the code presented in the Appendix  uses the proposed estimator to numerically approximate the marginal likelihood of this model with data simulated from a normal distribution. For example, running the simulation with $N=1000$ points sampled from the posterior, we find a numerical estimate $\log \widehat{f_X(x)} \approx -66.67655$ in good agreement with the theoretical value $\log f_X(x)\approx -66.67619$.

\section*{Appendix: computer code}

\bigskip

\linespread{1}

{\tiny
\begin{verbatim}
library(KernSmooth)

set.seed(1702)

si <- 3; n <- 25

x <- rnorm(n, mean = -1, sd = si)
m <- mean(x)
s <- ((n - 1) / n) * sd(x)

th0 <- 0; si0 <- 10

N <- 1000
th <- rnorm(N, mean = (n * si0^2 * m + si^2 * th0) / (n * si0^2 + si^2), 
                 sd = sqrt(si^2 * si0^2 / (n * si0^2 + si^2)))

post <- bkde(x = th)

log_f <- function(th) (-n / 2) * log(2 * pi) - n * log(si) +
                      (-n / (2 * si^2)) * ((th - m)^2 + s^2)

estimate <- mean(exp(log_f(th) + 
                 log(dnorm(th, mean = th0, sd = si0)) - log(approx(post$x, post$y, th)$y)))

f_X <- ((2 * pi)^(-n / 2) * si^(-n + 1) / sqrt(n * si0^2 + si^2)) * 
       exp((-1/2) * (n * s^2 / si^2 
                     - (n * si0^2 * m + si^2 * th0)^2 / (si^2 * si0^2 * (n * si0^2 + si^2)) 
                     + (n * si0^2 * m^2 + si^2 * th0^2) / (si^2 * si0^2)))

cat("Theoretical: ", log(f_X), "\n")
cat("Estimate:    ", log(estimate), "\n")
\end{verbatim}}

\bigskip\bigskip

\end{document}